\documentclass[onefignum,onetabnum]{siamonline171218}

\usepackage{tabularx} 	 
\usepackage{xspace}  	

\usepackage{lipsum}
\usepackage{amsfonts}
\usepackage{graphicx}
\usepackage{epstopdf}
\usepackage{algorithmic}
\usepackage{subcaption}
\usepackage{pgfplots} 
\pgfplotsset{compat=newest} 
\pgfplotsset{plot coordinates/math parser=false} 
\newlength\figureheight 
\newlength\figurewidth 
\usepackage{enumitem}
\setlist[enumerate]{leftmargin=.5in}
\setlist[itemize]{leftmargin=.5in}

\newsiamremark{remark}{Remark}
\newsiamremark{hypothesis}{Hypothesis}
\crefname{hypothesis}{Hypothesis}{Hypotheses}
\newsiamthm{claim}{Claim}

\headers{Gramians, EF and BT for Linear Systems with Quadratic Outputs}{P. Benner, P. Goyal and I. Pontes Duff}

\title{Gramians, Energy Functionals and Balanced Truncation for Linear Dynamical Systems with Quadratic Outputs\thanks{Submitted to the editors \today.
	}
}

\author{Peter Benner\thanks{Max Planck Institute for Dynamics of Complex Technical Systems, Sandtorstra\ss e 1, 39106 Magdeburg, Germany
		(\email{benner@mpi-magdeburg.mpg.de}).}
	\and Pawan Goyal\thanks{Corresponding author. Max Planck Institute for Dynamics of Complex Technical Systems, Sandtorstra\ss e 1, 39106 Magdeburg, Germany
		(\email{goyalp@mpi-magdeburg.mpg.de}).}
	\and Igor Pontes Duff\thanks{Max Planck Institute for Dynamics of Complex Technical Systems, Sandtorstra\ss e 1, 39106 Magdeburg, Germany
		(\email{pontes@mpi-magdeburg.mpg.de}).}
}

\usepackage{amsopn}
\DeclareMathOperator{\diag}{diag}

\usepackage{soul}

\newcommand*\diff{\mathop{}\!\mathrm{d}}
\newcommand{\LDSQO}{{\texttt{LD\_QO}}}
\newcommand{\LDS}{{\texttt{LD}}}
\usepackage[hardware,software]{mymacros}
\usepackage{amsmath}
\usepackage{breqn}
\usepackage{tikz,pgfplots} 
\usepackage{empheq}
\newlength\wex
\newlength\hex
\begin{document}

\maketitle

\begin{abstract}
  Model order reduction is a technique that is used to construct low-order approximations of large-scale dynamical systems. 
  In this paper, we investigate a balancing based model order reduction method for dynamical systems with a linear dynamical equation and a quadratic output function.   To this aim, we propose a new algebraic observability Gramian for the system based on Hilbert space adjoint theory. We then show the proposed Gramians satisfy a particular type of generalized Lyapunov equations and we investigate  their connections to energy functionals, namely, the controllability and observability. This allows us to find the states that are hard to control and hard to observe via an appropriate balancing transformation. Truncation of such states yields reduced-order systems. Finally, based on  $\mathcal H_2$ energy considerations, we, furthermore, derive error bounds, depending on the neglected singular values. The efficiency of the proposed method is demonstrated by means of two semi-discretized partial differential equations and is compared with the existing model reduction techniques in the literature.
\end{abstract}

\begin{keywords}
  model order reduction, dynamical systems, balanced truncation, controllability and observability, Gramians, energy functionals.
\end{keywords}

\begin{AMS}
  15A69, 34C20, 41A05, 49M05, 93A15, 93C10, 93C15.
\end{AMS}

\section{Introduction}
Dynamical models of real-world problems are often constructed, for example, with the purpose of simulations, predictions, optimization, and control. Such models are usually governed by partial differential equations (PDEs). More often than not, a spatial discretization of PDEs is necessary to perform engineering studies. To capture the important dynamics of a complex phenomenon, we require a fine spatial discretization, thus leading to a large number of equations. This imposes a huge computation burden. To overcome this issue,  one can use model order reduction, aiming at constructing lower-dimensional models which capture the important dynamical behaviors of the original large-scale dynamical system.  

A particular class of dynamical systems, which often occurs in modeling, is linear dynamical (\LDS) systems. Model order reduction (MOR) of such systems has been extensively studied in the literature and has been successfully applied for various real-world problems, see e.g., \cite{morAnt05,morBenCOetal17,morBenMS05}. In this paper, we consider a variant of \LDS~systems, as follows:
\begin{subequations}
	\begin{align*}
	\dot{x}(t) &=Ax(t) + Bu(t), \quad x(0) = 0,  \\
	y(t) &={x(t)}^{T}Mx(t),
	\end{align*}
\end{subequations}
where  $A \in \mathbb{R}^{n \times n}$, $B \in \mathbb{R}^{n \times m}$ and $M \in \mathbb{R}^{n \times n}$. As can be seen in the above equation,  the state  equation is linear, but  the output equation takes a quadratic form of the state as shown. We refer to \eqref{eqn:LDS_QO_SystemEquations}  as an \LDSQO~system. These kinds of systems, particularly, appears when one's interests lie in observing, e.g., the variance or deviation of the state variables from a reference point. This happens to be the case in random vibration analysis \cite{lutes2004random} and problems where response quantities related to energy or power are considered. 

MOR techniques  for \LDSQO~systems have been investigated in the past. For instance, the authors in \cite{Beeumen.2010, Beeumen.2012} have proposed to rewrite an \LDSQO~system as an \LDS~system. It is then followed by reducing by well-known techniques for \LDS~systems such as balanced truncation and interpolation-based methods. Furthermore, very recently, the authors in \cite{Pulch.2017} proposed an alternative approach, where the \LDSQO~system is written as a quadratic-bilinear (\texttt{QB}) system. Subsequently, the \texttt{QB} system is reduced by tailoring the balanced truncation approach, proposed in \cite{morBenG17}. However, in the formal approach, we do not directly utilize the quadratic structure of the output equation, and the latter approach is not only numerically expensive, but it also fails to keep the structure of the original system into reduced systems. 

In this paper, we study a balanced truncation method for \LDSQO~systems by proposing a novel pair of Gramians for the system. We, furthermore, investigate energy functionals, namely controllability, and observability, which form a ground for a balancing based procedure. We also  characterize the controllability and observability of the system based on the Gramians.  This allows us to construct reduced-order systems, removing less important subspaces for the dynamics. In this procedure, it is not required to rewrite the system as an \LDS~or \texttt{QB} system, and it, inherently, preserves the \LDSQO~structure in a reduced-order system.

A precise structure of the remaining paper is as follows. In the subsequent section, we provide the problem description and briefly revise the state-of-the-art. In \Cref{sec:SPBT_Description}, we present a novel pair of Gramians, the so-called controllability and observability Gramians for \LDSQO~systems. Based on these Gramians, we study  controllability and observability energy functionals, allowing us to determine  the states that are hard to reach as well as hard to observe. Consequently, we propose a balancing method to construct a good quality of reduced-order systems. In \Cref{sec:error_bounds},  we derive error bounds for the approximation error between the original and the reduced-order systems. These errors bounds rely on the definition of the $\cH_2$ norm for such class of systems. Precisely, our contributions in the aspect are twofold. Firstly, we derive a general a posteriori expression. Secondly, we show that the error bounds are directly associated with the neglected singular values. In \Cref{sec:experiments}, we test the accuracy of the proposed method and compare with the techniques proposed in \cite{Pulch.2017,Beeumen.2010}.  Finally, in \Cref{sec:conclusion}, we conclude the paper by listing our main contributions and provide future directions.


\section{Problem Formulation and Background Work}\label{sec:ProblemDefinition}
In this section, we discuss the MOR problem and its related work in the literature. We begin by recalling the classical balanced truncation method for \LDS~systems.

\subsection{Balanced Truncation of LD systems }\label{sec:LDS_Definition}
Let us consider an \LDS~system of the form:
\begin{subequations}\label{eqn:LDS_SystemEquations}
	\begin{align}
	\dot{x}(t) &= Ax(t) + Bu(t),\quad x(0) = 0, \\
	y(t) &= Cx(t),
	\end{align}
\end{subequations}
with $A \in \mathbb{R}^{n \times n}$, $B \in \mathbb{R}^{n \times m}$ and $C \in \mathbb{R}^{p \times n}$; $x(t) \in \Rn$, $u(t) \in \Rm$ and $y(t) \in \Rp$ are the state, input and output vectors, respectively; $n,m,$ and $p$ denote the state dimension or the order of the system, the number of inputs, and the number of outputs, respectively.  The main purpose of MOR is to construct a low-dimensional system, precisely an $r$-dimensional system with $r \ll n$, approximating the behavior of the system \eqref{eqn:LDS_SystemEquations}.  To construct a reduced-order system, we employ the Petrov-Galerkin framework. It consists in finding  two projection matrices $V, W \in \R^{n \times r}$ such that $W^T V = I_r$. This allows us to construct a reduced-order system as follows:
\begin{subequations}\label{eqn:LDS_SystemEquations_Red}
	\begin{align}
	\dot{\hx}(t) &= \hA\hx(t) + \hB u(t),\quad \hx(0) = 0, \\
	\hy(t) &= \hC\hx(t),
	\end{align}
\end{subequations}
where $\hA = W^TAV  \in \mathbb{R}^{r \times r}$, $\hB = W^TB \in \mathbb{R}^{r \times m}$ and $\hC = CV \in \mathbb{R}^{p \times r}$. Furthermore, the matrices $V$ and $W$ should be designed in a way that the reduced-order system meets desired goals, e.g., $y \approx \hat{y}$, meaning that the original and reduced-order systems should provide very similar outputs when excited by the same input signal.

There exist several MOR techniques for \LDS~systems and  we refer the reader to the books \cite{morAnt05, morBenCOetal17}  for more details. In this work, we focus on balanced truncation (BT), which was introduced in the control systems literature in \cite{morMoo81}. It mainly relies  on the controllability and observability energy functionals, see e.g., \cite{Adjfujimoto2002}. The controllability energy functional is defined as the minimal amount of energy required  to steer the system to zero from a given state.  On the other hand, the observability energy functional can be defined as the output energy generated by a non-zero initial condition.
For \LDS~systems, these functionals can be given by the functions using the controllability and observability Gramians, respectively denoted by $P$ and $Q$, see e.g.,~\cite{morAnt05}. The Gramians satisfy the following Lyapunov equations:
\begin{subequations}\label{eqn:LDS_LyapunovEquations}
	\begin{align}
	AP + PA^T + BB^T &= 0, \label{eqn:LDS_LyapunovEquationsP}\\
	A^T Q + QA + C^T C &= 0.\label{eqn:LDS_LyapunovEquationsQ} 
	\end{align}
\end{subequations}
The main principle of BT lies in determining the states that are simultaneously hard to reach and hard to observe, in other words, these states require a lot of energy to steer from zero, as well as,  generate very little output energy. To identify such states, we make use of the balancing tool based on the Gramians, leading to a reduced-order system on truncation.  Furthermore,  it preserves stability and provides guaranteed error bounds, see e.g., \cite{morAnt05}. In \Cref{alg:LTBT1_Algorithm}, we  sketch the square-root balanced truncation algorithm, enabling us to determine a reduced-order system.

It is worth mentioning that the most costly step in the algorithm is to compute the Gramians as the solutions of \cref{eqn:LDS_LyapunovEquations}. However, the solutions of the equations  exhibit low-rank phenomena, i.e., there exists $Z_P\in \R^{n \times l}$, with $l \ll n$, such that $P \approx Z_PZ_P^T$. Thus, we make use of low-rank solvers of Lyapunov equations, which are not only numerically efficient but also yield directly the Gramians in Cholesky factors.  In the past decades, several advance algorithms have been proposed, allowing us to compute a solution of a Lyapunov equation of a few thousands in a low-rank form on a moderate machine. We refer the reader to the review papers~\cite{BenLP08, BenS13,Sim16a} to get an overview of the existing methods. 
\begin{algorithm}[tb!]
	\caption{Balanced truncation method for \LDS~systems.}
	\label{alg:LTBT1_Algorithm}
	\begin{algorithmic}[1]
		\REQUIRE{Matrices $(A,B, C)$  and the order of a reduced-order system $r$.}
		\ENSURE{Reduced matrices $(\hA, \hB, \hC)$.}
		\STATE Compute  low-rank factors of Gramians  $P \approx Z_PZ_P^T$, $Q \approx Z_QZ_Q^T$, where $P$ and $Q$ solve for \eqref{eqn:LDS_LyapunovEquationsP} and \eqref{eqn:LDS_LyapunovEquationsQ}, respectively.
		\STATE Compute the SVD of $Z_P^TZ_Q$, and partition as follows: \[ Z_P^TZ_Q  = \begin{bmatrix}
		U_1 & U_2
		\end{bmatrix}\diag \left(\Sigma_1, \Sigma_2\right)\begin{bmatrix}
		V_1 & V_2
		\end{bmatrix}^T,\quad \text{with}~ \Sigma_1 \in \R^{r\times r}. \vspace{-0.5cm}\] 
		\STATE Construct the projection matrices $ V = Z_{P}U_1\Sigma_1^{-\frac{1}{2}}$ and $W = Z_{Q}V_1\Sigma_1^{-\frac{1}{2}}.$  
		\STATE Construct $\hat{A} = W^T {A} V$,\quad  $\hat{B} = W^T {B}$,\quad  $\hat{C} = {C} V$.
		\RETURN$ \hA, \hB,$ and $\hC$.
	\end{algorithmic}
\end{algorithm}

\subsection{Problem formulation for LDSQO systems}
In this paper, our focus rather lies on linear dynamical systems with quadratic output function (\LDSQO~systems) which is of the form:
\begin{equation}\label{eqn:LDS_QO_SystemEquations} 
H   :=
\begin{cases}
\dot{x}(t) =Ax(t) + Bu(t), \quad x(0) = 0, 	 \\
y(t) ={x(t)}^{T}Mx(t), 			
\end{cases}
\end{equation}
where $A \in \mathbb{R}^{n \times n}$, $B \in \mathbb{R}^{n \times m}$ and $M \in \mathbb{R}^{n \times n}$. For ease, the original \LDSQO~system~\eqref{eqn:LDS_QO_SystemEquations} is denoted by $H = (A, B, M)$.    We assume that the matrix $A$ is Hurwitz; hence, the \texttt{LDS_QO} system is asymptotically stable.  Furthermore, without loss of generality, we assume that the matrix $M$ is symmetric, i.e., $M = M^T$. In case the matrix $M$ is not symmetric, we can always construct the symmetric matrix  $M_s := \left(M+M^T\right)/2$, which ensures $x^T(t)Mx(t) = x(t)^TM_s x(t)$. 
Note that an \LDSQO~system \eqref{eqn:LDS_QO_SystemEquations} has an input-output nonlinear mapping, even if the dynamical equation is linear. Our aim is to find two projection matrices $V, W \in \Rnr$, with $W^TV = I_r$, allowing to construct a reduced-order system $\hH := (\hA, \hB, \hM)$ as follows:
\begin{equation}\label{eqn:LDS_QO_SystemEquations_Red} 
\hH  : =
\begin{cases}
\dot{\hx}(t) =\hA\hx(t) + \hB u(t),\quad \hat x(0) = 0, \\
\hy(t) ={\hx(t)}^{T}\hM\hx(t),
\end{cases}
\end{equation}
where $\hA = W^TAV \in \mathbb{R}^{r \times r}$, $\hB = W^TB \in \mathbb{R}^{r \times m}$ and $\hM = V^TMV\in \mathbb{R}^{r \times r}$ while ensuring the desired properties. 
\subsection{State of the Art}\label{sec:StateOfTheArt}
In this subsection, we briefly discuss two existing MOR methods for \LDSQO~systems, see \cite{Pulch.2017, Beeumen.2010}. 
\subsubsection{Linear transformation and balanced truncation}\label{sec:LTBT_Description}
A MOR method for \LDSQO~systems was proposed in~\cite{Beeumen.2010}. The primary concept of the method is to equivalently rewrite an \LDSQO~system as an \LDS~system. This allows us to employ BT  for the \LDS~system. As discussed in \cite{Beeumen.2010}, for $M\succeq 0$, we can rewrite the system \eqref{eqn:LDS_QO_SystemEquations} as follows:
\begin{subequations}\label{eqn:LTBT_SystemEquations}
	\begin{align}
	\dot{{x}}(t) &= {A} {x}(t) + Bu(t),\quad x(0) = 0, \\
	y_{\texttt T}(t) &=  C_{\texttt T}x(t), 
	\end{align}
\end{subequations}
where the matrices $A, B$ and the state $x(t)$ are the same as in \eqref{eqn:LDS_QO_SystemEquations}, but the matrix $C_{\texttt T} \in \R^{q \times n}$ is such that $C_{\texttt T}^TC_{\texttt T} = M$, where $q = \rank{M}$. However, the outputs of the systems \eqref{eqn:LDS_SystemEquations} and \eqref{eqn:LDS_QO_SystemEquations} are related as $y(t) = \|y_{\texttt T}(t)\|_2^2$. This transformation allows us to employ the BT method for the \LDS~system. For completeness, we sketch the steps to determined a reduced-order system using the procedure given in \cref{alg:LTBT_Algorithm}.  One of drawbacks of the method is that it does not make use of the quadratic-form of the output equation. Furthermore, later in the paper, we will discuss that the approach might be very expensive when $q$ is large or when $M \nsucceq 0$.

\begin{algorithm}[!tb]
	\caption{BT for \LDSQO~systems (rewriting as an \LDS~system from~\cite{Beeumen.2010})}
	\label{alg:LTBT_Algorithm}
	\begin{algorithmic}[1]
		\REQUIRE{The original system's  m3atrices $(A,B, M)$  and  the order of the reduced-order system~$r$.}
		\ENSURE{Reduced-order system's matrices $(\hA, \hB, \hC)$.}
		\STATE{Compute the Cholesky factorization of the matrix $M = C_{\texttt T}^TC_{\texttt T}$.}
		\STATE{For the \LDS~system given by $A, B, C_{\texttt T}$, apply \Cref{alg:LTBT1_Algorithm} to obtain $\hA, \hB, \hC.$ }
		\RETURN$ \hA, \hB,$ and $\hC$.
	\end{algorithmic}
\end{algorithm}

\begin{remark}
	\Cref{alg:LTBT_Algorithm} is shown when the matrix $M$ is symmetric positive semi-definite. However, the authors in \cite{Beeumen.2010} also discusses the case general case where $M$ is not a positive semi-definite. We refer to the reference for details.
\end{remark}


\subsubsection{Quadratic-bilinear transformation and balanced truncation}\label{sec:QBTBT_Description}
Recently, a novel approach to construct reduced-order systems for \LDSQO ~systems has been proposed in \cite{Pulch.2017}. This approach consists of, firstly, converting the original \LDSQO~system as a quadratic-bilinear (\texttt{QB}) system by taking the derivative of the output equation of the original \LDSQO~system and followed by augmenting it into the state variable $x(t)$. It results into the \texttt{QB} system as follows:
\begin{subequations}
	\label{eqn:QBTBT_SystemEquations}
	\begin{align}
	\dot{x}_{\texttt{qb}}(t) &={A}_{\texttt{qb}}{x}_{\texttt{qb}}(t) + {B}_{\texttt{qb}}{u}(t) 
	+ {H}_{\texttt{qb}}({x}_{\texttt{qb}} \otimes {x}_{\texttt{qb}}) + 
	\sum_{j = 1}^{m} {{u}_j} {N}_{\texttt{qb}} {{x}_{\texttt{qb}}},\quad x_{\texttt{qb}}(0) = 0,\\
	y(t) &= C_{\texttt{qb}}x_{\texttt{qb}}(t),
	\end{align}
\end{subequations}
where 
\begin{align*}
{x}_{\texttt{qb}} &= \begin{bmatrix}x\\ y\end{bmatrix}, & A_{\texttt{qb}} &= \begin{bmatrix}A&0\\0&0\end{bmatrix},&
{H}_{\texttt{qb}} &= \begin{bmatrix}0&0&0&0&\cdots&0&0&0&0\\s_{1}^{T}&0&s_{2}^{T}&0&\cdots&s_{n}^{T}&0&0&0\end{bmatrix}\\
{B}_{\texttt{qb}} &= \begin{bmatrix}B\\0\end{bmatrix}, & C_{\texttt{qb}} &= \begin{bmatrix}1&0\end{bmatrix}, & {N}^{(j)}_{\texttt{qb}} &= \begin{bmatrix}0&0\\2 b{_	j}^T M &0\end{bmatrix}, ~j = 1,\ldots,m
\end{align*}
in which  $x(t), u(t),y(t), A,B$, and $C$ are as defined in \eqref{eqn:LDS_QO_SystemEquations}; $b_j \in \Rn$ is the $j$-th column vectors of $B$; $s_j \in \mathbb{R}^{n} $ is the $j$-th column of the the matrix $S$, which is defined as follows:
\begin{equation}\label{eqn:QBTBT_SymmetricMatrixEquation}
S =A^T M + M^T A.
\end{equation}
The system \eqref{eqn:QBTBT_SystemEquations} has $n+1$ states, $m$ inputs, and $1$ output. Note that $A_{\texttt{qb}}\in \mathbb{R}^{(n+1) \times (n+1)}$, $B_{\texttt{qb}} \in \mathbb{R}^{(n+1) \times m}$, $N^{(j)}_{\texttt{qb}} \in \mathbb{R}^{m \times (n+1)}$ and $H_{\texttt{qb}} \in \mathbb{R}^{(n+1) \times (n+1)^2}$.

Once we have the equivalent \texttt{QB} system, we can employ recently developed MOR schemes for \texttt{QB} systems, see e.g., \cite{morBenB15,morBenG17,morBenGG18}. Focusing on a BT method,  controllability Gramian $P_{\texttt{qb}}$ and the observability Gramian $Q_{\texttt{qb}}$ for a \texttt{QB} system was proposed in \cite{morBenG17}, where it has been shown that these Gramians satisfy  the following quadratic-type Lyapunov equations:
\begin{subequations}\label{eqn:QBTBT_LyapunovEquations}
	\begin{align}
	{A}_{\texttt{qb}} P_{\texttt{qb}} + P_{\texttt{qb}} {A}_{\texttt{qb}}^T + {B}_{\texttt{qb}}{B}_{\texttt{qb}}^T 
	+ {H}_{\texttt{qb}}\left(P_{\texttt{qb}} \otimes P_{\texttt{qb}}\right){H}_{\texttt{qb}}^T  + 
	\sum_{j = 1}^{m} N^{(j)}_{\texttt{qb}} P_{\texttt{qb}} \left(N^{(j)}_{\texttt{qb}}\right)^T&= 0,\label{eqn:QBTBT_LyapunovEquations_Reach}\\
	{A}_{\texttt{qb}}^T Q_{\texttt{qb}} + Q_{\texttt{qb}} {A}_{\texttt{qb}} + {C}_{\texttt{qb}}^T{C}_{\texttt{qb}}
	+ {H}_{\texttt{qb}}^{(2)}\left(P_{\texttt{qb}} \otimes Q_{\texttt{qb}}\right)\left({H}_{\texttt{qb}}^{(2)}\right)^T  + 
	\sum_{j = 1}^{m} \left({N^{(j)}_{\texttt{qb}}}\right)^T  Q_{\texttt{qb}} N^{(j)}_{\texttt{qb}}&= 0. \label{eqn:QBTBT_LyapunovEquations_Obser}
	\end{align}
\end{subequations}

Having had the Gramians for \texttt{QB} systems, one can obtain a reduced-order system using the classical square-root method, as shown in \cite{morBenG17}. However, we would like to list an additional challenge while solving \eqref{eqn:QBTBT_LyapunovEquations}, that is the matrix $A_{\texttt{qb}}$ is a singular matrix, i.e., it contains zero eigenvalues. This issue has been addressed in \cite{Pulch.2017}, where the authors have discussed in details on how to solve \eqref{eqn:QBTBT_LyapunovEquations} efficiently. We sketch the proposed methodology in  \Cref{alg:QBTBT_Algorithm}.

Both methodologies presented in \cite{Pulch.2017,Beeumen.2010} rely on rewriting the original \LDSQO~system into a new form and then apply an appropriate MOR method to construct ROMs.
In the next section, we proposed Gramians and a MOR  method which relies on the original structure of the \LDSQO~system. 

\begin{algorithm}[!tb]
	\caption{QB transformation and BT for \LDSQO ~systems.} 
	\label{alg:QBTBT_Algorithm}
	\begin{algorithmic}[1]
		\REQUIRE{The original system's matrices $(A,B, M)$  and the order of a reduced-order system $r$.}
		\ENSURE{Reduced QB system's matrices $\hat{A}_{\texttt{qb}}, \hat{B}_{\texttt{qb}}$, $\hC_{\texttt{qb}}$, $\hH_{\texttt{qb}}$, and $\hat{N}_{\texttt{qb}}^{(j)}$.}
		\STATE{Determine QB system matrices, i.e., $A_{\texttt{qb}}$, $B_{\texttt{qb}}$, $C_{\texttt{qb}}$, $H_{\texttt{qb}}$ and $N^{(j)}_{\texttt{qb}}$, as shown in \eqref{eqn:QBTBT_SystemEquations}.}
		\STATE{Compute  low-rank factors of Gramians 
			$P_{\texttt{qb}} \approx Z_{P_{\texttt{qb}}} Z_{P_{\texttt{qb}}}^T$, $Q_{\texttt{qb}} \approx Z_{Q_{\texttt{qb}}} Z_{Q_{\texttt{qb}}}^T$, 
			where $P$ and $Q$ solve \eqref{eqn:QBTBT_LyapunovEquations_Reach} and \eqref{eqn:QBTBT_LyapunovEquations_Obser}, respectively.}\vspace{-0.0cm}
		\STATE{Compute the SVD of $Z_{P_{\texttt{qb}}}^TZ_{Q_{\texttt{qb}}}$, and decompose as: \[ Z_{P_{\texttt{qb}}}^TZ_{Q_{\texttt{qb}}} = \begin{bmatrix}
			U_1 & U_2
			\end{bmatrix}\diag{\left(\Sigma_1, \Sigma_2\right)}\begin{bmatrix}
			V_1 & V_2
			\end{bmatrix}^T,\quad \text{with}~ \Sigma_1 \in \R^{r\times r}. \vspace{-0.5cm}\] }
		\STATE{ Construct the projection matrices $ V = Z_{P_{\texttt{qb}}}U_1\Sigma_1^{-\frac{1}{2}}$ and $W = Z_{Q_{\texttt{qb}}}V_1\Sigma_1^{-\frac{1}{2}}.$  
		}
		\STATE Construct reduced-order matrices:
		\begin{align*}
		\hat{A}_{\texttt{qb}} &= W^T {A}_{\texttt{qb}} V,&   \hat{B} &= W^T B_{\texttt{qb}}, &\hat{N}^{(j)}_{\texttt{qb}} &= W^T N^{(j)}_{\texttt{qb}} V, \quad j = 1,\ldots,m, \\
		\quad\hat{C}_{\texttt{qb}}&=C_{\texttt{qb}}V,&\hat{H}_{\texttt{qb}} &= W^T {H}_{\texttt{qb}}\left(V\otimes V\right).
		\end{align*}
		\RETURN{Reduced \texttt{QB}  system's matrices $\hat{A}_{\texttt{qb}}, \hat{B}_{\texttt{qb}}$, $\hC_{\texttt{qb}}$, $\hH_{\texttt{qb}}$, and $\hat{N}_{\texttt{qb}}^{(j)}$.}
	\end{algorithmic}
\end{algorithm}

\section{Balanced truncation method for LDS\_QO systems} \label{sec:SPBT_Description}
This section contains the main theoretical contributions of the paper. We  focus on deriving a new pair of Gramians for the \LDSQO~system. In particular, we propose a tailored observability Gramian  for the considered systems by means of the adjoint systems for nonlinear systems \cite{Adjfujimoto2002}. We begin by discussing the controllability Gramian for \LDSQO~systems. 
\subsection{Controllability  Gramian and controllability energy functional} Since the differential equation of a \texttt{LDS_QO} system is a first-order linear time-invariant equation, it is well-known that the controllability  Gramian $P$ is defined as follows:
\begin{equation}
P := \int_0^\infty e^{A\tau}BB^Te^{A^T\tau}d\tau.
\end{equation}
Moreover, if $A$ is Hurwitz, i.e., $\sigma(A) \subset \C^{-}$, then the controllability Gramian $P$ satisfies the following Lyapunov equation:
\begin{equation}\label{eqn:LDSQO_LyapunovEquations_Reach}
AP + PA^T + BB^T = 0.
\end{equation}
Moreover, the controllability energy functional $E_c(x_0)$ is defined as minimum input energy required to steer the state from a non-zero initial condition to zero, i.e.,
\begin{equation*}
E_c(x_0) = \min_{\begin{array}{c}\scriptstyle x(-\infty) = x_0, \\[-4pt]
	\scriptstyle x(0) = 0
	\end{array}} \|u\|^2_{L_2}.
\end{equation*}
Furthermore, the controllability energy functionals can be given in terms of the controllability Gramian as follows:
\begin{equation*}
E_c(x_0) = \dfrac{1}{2}x_0^TP^{-1}x_0,
\end{equation*}
assuming $P>0$.  The controllability energy functional relation shows us that the state components, corresponding to the smaller singular values of the Gramian $P$ are hard to reach. We refer, e.g., to \cite{morAnt05} for more details. 
\subsection{Observability Gramian}
A major difference between a classical first-order linear time-invariant  and \LDSQO~system is the output equation, i.e., the formal system has the output equation as $Cx(t)$, whereas the output equation of the latter system  takes a quadratic form, given by $x(t)^TMx(t)$. Hence, we expect to have a different observability Gramian, which somehow relates the output energy functional of the \LDSQO~system. For this, we make use of the adjoint theory for nonlinear systems, developed in \cite{Adjfujimoto2002}. Following the discussion given in the listed reference, we can write down the state-space realization of the nonlinear Hilbert adjoint operator of an \LDSQO~system as follows:
\begin{subequations}
	\begin{align}
	&\dot{x}(t) = Ax(t) +Bu(t),~\,~x(0) = 0,  \label{eq:StateEq}\\ 
	&\dot{z}(t) = -Az(t) -Mx(t)u_d(t),~\,~z(\infty) = 0, \label{eq:AdjEqState} \\
	&y_d(t) = B^Tz(t),
	\end{align}
\end{subequations}
where $z(t) \in \R^n,$ $u_d(t)\in \R$ and $y_d(t) \in \R^m$ are, respectively, the dual state, the dual input and the dual output. Based on this Hilbert adjoint operator, we construct the transfer map between the adjoint input and the adjoint state. To that end, let us integrate the adjoint equation \eqref{eq:AdjEqState} backwards as follows:
\[ z(t) = \int_{\infty}^t e^{-A^T(t-\sigma)}Mx(\sigma)u_d(\sigma)\diff\sigma.  \]
After a suitable change of variables, we obtain
\[ z(t) = 
\int_{\infty}^0 e^{A^T\sigma_1}Mx(t+\sigma_1)u_d(t+\sigma_1)\diff\sigma_1.\]
Moreover, from \eqref{eq:StateEq}, we have
\[x(\sigma_1+t) = \int_0^{\sigma_1+t} e^{A\sigma_2}Bu(t-\sigma_2)\diff\sigma_2. \]
By injecting the above expression in the former equation, we obtain
\begin{align}\label{eq:dualstate}
z(t) = 
&\int_{\infty}^0\int_0^{\sigma_1+t} e^{A^T\sigma_1}M e^{A\sigma_2}Bu(t-\sigma_2)u_d(t+\sigma_1)\diff\sigma_2\diff\sigma_1.
\end{align}
This allows us to define the observability Gramian of the \LDSQO~system as follows:
\begin{equation}\label{eq:observabilityGramian}
Q = 
\int_0^\infty\int_0^\infty e^{A^T\sigma_1}M e^{A\sigma_2}B\left(e^{A^T\sigma_1}M e^{A\sigma_2}B\right)^Td\sigma_1d\sigma_2.
\end{equation}
In what follows, we present the matrix equation, solving for the observability Gramian $Q$, defined in \eqref{eq:observabilityGramian}. 
\begin{lemma}
	Let the observability Gramian $Q$ be defined as in \eqref{eq:observabilityGramian}. Assuming the matrix $A$ in an {\LDSQO} system is Hurwitz, the observability Gramian $Q$ is a unique solution to the following Lyapunov equation:
	\begin{equation}\label{eqn:LDSQO_LyapunovEquations_Obser}
	A^TQ + AQ + MPM =0,
	\end{equation}	
	where $P$ is the controllability Gramian, satisfying 
	\begin{equation}
	AP + PA^T + BB^T =0.
	\end{equation}	
\end{lemma}
\begin{proof}
	From \eqref{eq:observabilityGramian}, we know the observability Gramian satisfies:
	\begin{align}
	Q &= \int_0^\infty\int_0^\infty e^{A^T\sigma_1}Me^{A\sigma_2}B\left(e^{A^T\sigma_1}Me^{A\sigma_2}B\right)^Td\sigma_1d\sigma_2\nonumber\\
	&= \int_0^\infty e^{A^T\sigma_1}M\left(\int_0^\infty e^{A\sigma_2}BB^T e^{A^T\sigma_2}\right)Me^{A\sigma_1} d\sigma_1d\sigma_2.\label{eq:QwithP1}
	\end{align}
	We know from, e.g., \cite{morAnt05}, that 
	\begin{equation}\label{eq:eqnforP}
	\int_0^\infty e^{A\sigma_2}BB^T e^{A^T\sigma_2}d\sigma_2 = P,
	\end{equation}
	where $P$ is the controllability Gramian of \LDSQO~systems. Substituting the above relation into \eqref{eq:QwithP1} yields 
	\begin{align}
	Q &=  \int_0^\infty e^{A^T\sigma_1}M PMe^{A\sigma_1} d\sigma_1.\label{eq:QwithP2}
	\end{align}
	Using the same arguments as used for \eqref{eq:eqnforP}, for $\sigma(A) \in \C_{-}$, it can be readily shown that  the observability Gramian $Q$ satisfies \eqref{eqn:LDSQO_LyapunovEquations_Obser}.
\end{proof}
Next, we investigate a relation between the observability Gramian and energy functionals. The observability energy functional $E_o(x_0)$ is defined as the output energy produced by  the nonzero initial condition $x_0$, i.e., 
\[ E_o(x_0) = \int_0^{\infty} \|y(t)\|^2 dt, \]
where $y(t)$ is the output of a system. To that end, we establish a relation between the observability Gramian and observability energy functionals $E_o(x_0)$ in the following theorem.  
\begin{theorem}
	Let the controllability $P> 0$ and observability $Q$ be defined as \eqref{eqn:LDSQO_LyapunovEquations_Reach} and \eqref{eqn:LDSQO_LyapunovEquations_Obser}, respectively. Furthermore, let us assume that the state trajectory $x(t)$, generated from a non-zero initial condition $x_0$ with $u(t)  \equiv 0$, lies in $\cW_\delta$, where $\cW_\delta$ denotes the balls of radius $\delta$ centered around zero.  Then, the output energy functional can be bounded as follows:
	\begin{equation*}
	E_o(x_0) \leq x_0^T Qx_0.
	\end{equation*}
\end{theorem}
\begin{proof} 
	Using the definition of the observability energy functional, we have
	\begin{align*}
	E_o(x_0) &= \int_0^t \|y(t)\|_2^2dt\\
	&= \int_0^t x(\tau)^TM^Tx(\tau) x(\tau)^TMx(\tau)d\tau.
	\end{align*}
	Note that $x(\tau)$ can be given as $e^{A\tau}x_0$. Since the system is assumed to be controllable ($P>0$), we can write $x(\tau) = Lz(\tau)$ $\forall \tau$, where the matrix $L$ is the Cholesky factors of $P$, i.e., $L L^T = P$. Using all these relations, we obtain
	\begin{align*}
	E_o(x_0) &= \int_0^t  x_0e^{A^T\tau}M^T {e^{A\tau}x_0x_0^T e^{A^T\tau}}Me^{A\tau}x_0 d\tau \\
	&= \int_0^t  x_0e^{A^T\tau}M^T Lz(\tau)z(\tau)^TL^TMe^{A\tau}x_0 d\tau.
	\end{align*}
	Additionally, it can be easily shown that $Lz(t)z(t)^TL^T \leq LL^T$ for $\|z(t)\|^2_2 \leq 1$. Hence, if an initial condition $x_0$ is such that the generated state trajectory $x(\tau) \in \cW_\delta$, where $\cW_\delta$ is chosen such that every $x(\tau)$ can be written as $Lz(t)$, where $\|z(t)\|^2_2\leq 1$. Thus, we get
	\begin{align*}
	E_o(x_0) &= \int_0^t  x_0e^{A^T\tau}M^T Lz(\tau)z(\tau)^TL^TMe^{A\tau}x_0 d\tau\\
	&\leq  \int_0^\infty  x_0^Te^{A^T\tau}M^TLL^TMe^{A\tau}x_0 d\tau\\
	&\leq  \int_0^\infty  x_0^Te^{A^T\tau}M^TPMe^{A\tau}x_0 d\tau\\
	&\leq x_0^T\left( \int_0^\infty  e^{A^T\tau}M^TPMe^{A\tau} d\tau\right)x_0\\
	& \leq  \int_0^\infty  x_0^TQx_0.
	\end{align*}
	This concludes the proof.
\end{proof}

So far, we have proposed Gramians for {\LDSQO} systems and have shown how these Gramians relate to the energy functionals of the systems, under required conditions. However, in the following, we show that these Graimans, in a general case, encode controllability and observable subspaces information. 
\begin{theorem} \label{thm:ContObser_PQ}
	Let controllability Gramian ($P$) and observability Gramian ($Q$) solve \cref{eqn:LDSQO_LyapunovEquations_Reach} and \cref{eqn:LDSQO_LyapunovEquations_Obser}, respectively. Then, we have the following results:
	\begin{enumerate}
		\item[(a)] If the system is aimed at steering from zero to $x_0$, which belongs to $\ker{P}$, then $E_c(x_0) = \infty$; hence, it is unreachable.
		\item[(b)] If $P >0$ and the initial condition $x_0 \in \ker{Q}$, then $E_0(x_0) = 0$, thus making the state $x_0$ unobservable.
	\end{enumerate}
\end{theorem}
\begin{proof}
	(a) This result is very well-known in the literature; hence, for the brevity of the paper, we skip the details and refer the reader e.g., to \cite{morBenD11,morBenG17}, where authors have considered a more general case.
	
	(b) We know that the observability Gramian satisfies the following relation:
	\begin{equation}\label{eqn:observability1}
	A^TQ + QA + MPM = 0.
	\end{equation}
	Next, let us consider a vector $v \in \ker{Q}$ and multiply \eqref{eqn:observability1} from the left and right-hand sides by $v^T$ and $v$, respectively, yielding
	\begin{align*}
	v^TA^TQv + v^TQAv + v^TMPMv &= 0\\
	v^TMPMv &= 0.
	\end{align*}
	This implies $PMv = 0$. Furthermore, it can be noticed that $QAv = 0$. Next, we consider that $x(t) \in \ker{Q}$ at time $t$ and a vector $\tv \in \range{Q}$; hence, we have
	\begin{equation}
	\tv^T \dot {x}(t) = \tv^TAx(t) = 0.
	\end{equation}
	This means that if $x(t) \in \ker{Q}$, then $\dot {x}(t) \in \ker{Q}$. So, if the initial condition $x_0 \in \ker{Q}$, then $x(t) \in \ker{Q},\forall t \geq 0$. Furthermore, note that the system is assumed to be controllable; this means that $x(t) \in \range{P}$, i.e., $x(t) = P\tilde x(t), \forall t \geq 0$. Thus, the output $y(t)$ of the \LDSQO~system is given as $x(t)^TMx(t) = \tx(t)^TPMx(t) = 0$ since $x(t)$ also lies in $\ker{Q}$. As a result, the output energy functional is zero; hence, the initial state $x_0$ cannot be observed. 
\end{proof}
Having had all this discussion between the energy functionals and Gramians, it is clear that these Gramians allow us to determine the states which are hard to reach and hard to observe. In what follows, we propose a new BT algorithm for  \LDSQO~systems.
\subsection{New balanced truncation method for LDS\_QO systems}\label{subsec:SPBT}
The main idea of BT lies in furthermore neglecting the states which are both hard to reach and hard to observe states. 
In order to guarantee that hard to reach and hard to observe states are truncated simultaneously, we need to find a state transformation $T_{\mathcal{B}}$ such that  the \LDSQO~system is transformed into a balanced  realization. Thus, the controllability and observability Gramians of the transformed realization are the same and diagonal, i.e.,
\[P = Q = \Sigma = \diag(\sigma_1, \sigma_2, \dots, \sigma_n),   \]
where $\sigma_1\geq \sigma_2 \geq \dots \geq \sigma_n > 0$ and $\sigma_k$ are referred to as the singular values of \LDSQO~systems.  Such a  transformation exists whenever $P$ and $Q$ are positive definite matrices. Moreover, the small singular values $\sigma_k$ characterize the states that are hard to reach and hard to observe, which can then be truncated.  Next, we assume that the matrices of the balanced system $H = (A, B, M)$ are partitioned as
\begin{align}\label{eq:PartionedMatrices} A = \begin{bmatrix}
A_{11} & A_{12}  \\ A_{21} & A_{22}  
\end{bmatrix},\quad B = \begin{bmatrix}
B_{1}  \\  B_{2}  
\end{bmatrix}, \quad M = \begin{bmatrix}
M_{11} & M_{12}  \\ M_{12}^T & M_{22}  
\end{bmatrix}~\textnormal{and} \quad \Sigma = \begin{bmatrix}
\Sigma_1 & 0  \\ 0 & \Sigma_2 
\end{bmatrix},  \end{align}
where $\Sigma_1 = \diag(\sigma_1, \dots, \sigma_r)$ and $\Sigma_2 = \diag(\sigma_{r+1}, \dots, \sigma_n)$. Since the system $H$ is assumed to be balanced, we have
\begin{subequations}\label{eq:FullBalGram}
	\begin{align}
	A\Sigma +\Sigma A^T +BB^T = 0, \label{eq:FullReachGram}\\
	A^T\Sigma+\Sigma A +M\Sigma M = 0. \label{eq:FullObsevGram}
	\end{align}
\end{subequations} 

Subsequently, the reduced-order system can be easily obtained by considering the upper-left blocks, yielding  $\hH = (A_{11}, B_{1}, M_{11})$. By simple algebra, it can also be seen that  the reduced matrices satisfy the following equations:
\begin{subequations}
	\begin{align}
	A_{11}\Sigma_1 +\Sigma_1A_{11}^T +B_1B_1^T& = 0, \label{eq:RedRechGram_Part}  \\
	A_{11}^T\Sigma_1 +\Sigma_1A_{11} +M_{11}\Sigma_1 M_{11} +M_{12}\Sigma_2 M_{12}^T&=0, \label{eq:RedObservGram_Part}
	\end{align} 
\end{subequations}
thus allowing us to make the following observation.
\begin{remark} From \cref{eq:RedObservGram_Part}, one can conclude that the reduced-order system might not be balanced even if the original model is. A necessary condition for the reduced-order system to be a balanced one is $M_{12} = 0$, which in general is not true. 
	Nevertheless, since $M_{12}\Sigma_2M_{12}^T$ is a symmetric positive semi-definite matrix, the following matrix inequality holds
	\[ A_{11}^T\Sigma_1 + \Sigma_1A_{11} +M_{11}\Sigma_1M_{11} \leq 0.   \] 
	As a consequence, the  reduced-order system is balanced in the generalized sense, see \cite[Sec. 4.7]{dullerud2013course}.
\end{remark}

Given a system $H$, it is not necessary that the system is in a balanced form. One way to approximate it is to compute a balanced realization, which is followed by computing a reduced-order system as described above. However, analogous to the linear case, the balanced transformation is not required explicitly. Instead, one can construct two projection matrices $V$ and $W$ using the Cholesky factors of $P$ and $Q$ to determine directly a reduced-order system. This procedure is known as square-root BT, see \cite[Sec. 7.3]{morAnt05};
we sketch the steps to construct reduced-order systems for \LDSQO~systems based on the proposed Gramians in  \cref{alg:SPBT_Algorithm}.


\begin{algorithm}[tb]
	\caption{Novel BT method for \LDSQO~systems.}
	\label{alg:SPBT_Algorithm}
	\begin{algorithmic}[1]
		\REQUIRE{The original system's matrices $(A,B, M)$  and the order of the reduced-order system $r$.}
		\ENSURE{The reduced-order system's matrices $(\hA, \hB, \hM)$.}
		\STATE{Compute low factors of Gramians $P_{\texttt{}}$ and $Q_{\texttt{}}$, i.e., $P_{\texttt{}} \approx Z_{P_{\texttt{}}} Z_{P_{\texttt{}}}^T$ and $P_{\texttt{}} \approx Z_{Q_{\texttt{}}} Z_{Q_{\texttt{}}}^T$, where $P_{\texttt{}}$ and $Q_{\texttt{}}$ solve \eqref{eqn:LDSQO_LyapunovEquations_Reach} and \eqref{eqn:LDSQO_LyapunovEquations_Obser}, respectively.}
		\STATE Perform the SVD of $Z_P^TZ_Q$, and decompose as \[ Z_P^TZ_Q  = \begin{bmatrix}
		U_1 & U_2
		\end{bmatrix}\diag{\left(\Sigma_1, \Sigma_2\right)}\begin{bmatrix}
		V_1 & V_2
		\end{bmatrix}^T,\quad \text{with}~ \Sigma_1 \in \R^{r\times r}. \vspace{-0.5cm}\] 
		\STATE Construct the projection matrices $ V = Z_{P}U_1\Sigma_1^{-\frac{1}{2}}$ and $W = Z_{Q}V_1\Sigma_1^{-\frac{1}{2}}.$  
		\STATE Construct 	$\hat{A} = W^T A V$, \quad $\hat{B} = W^T A$,  \quad$\hat{M} = V^T M V$ . 
		\RETURN$ \hA, \hB,$ and $\hM$.
	\end{algorithmic}
\end{algorithm}

\subsection{Advantages of the proposed method}
Next, we note advantages of the proposed method over the existing BT methods for \LDSQO~systems, which are:
\begin{itemize}
	\item The methodology does not require any prior transformation of an \LDSQO~system into a classical linear system or a QB system. Hence, computational efforts converting it into an equivalent linear or \texttt{QB} system can be saved. 
	\item If an \LDSQO~system is written as a linear system, then for the observability Gramian, we need to solve $$QA^T+AQ = -C_T^TC_T,$$
	where $C_T^TC_T = M$ and $M\succeq 0$, whereas for \LDSQO~systems, we solve \eqref{eqn:LDSQO_LyapunovEquations_Obser} for the observability Gramian. So, note that the rank of the matrix ${M} P {M}$, where $P$ is the controllability Gramian, is always smaller than or equal to the rank of $C_T$, where $C_T^TC_T = M$ if $M\geq 0$. It is because $M PM$ can be seen as a projection of the controllable subspace onto the range of the matrix $M$.

	Broadly speaking, we know that a lower rank of the right-hand side of a Lyapunov equation can lead to a faster converge and also the solution is generally of a lower rank. Hence, solving \eqref{eqn:LDSQO_LyapunovEquations_Obser} might be computationally  efficient, thus constructing reduced-order systems. To illustrate this, we consider a $2$-dimensional dynamical system as follows:
	\begin{subequations}\label{eq:2by2example}
		\begin{align}
		\begin{bmatrix} \dot x_1(t) \\ \dot x_2(t) \end{bmatrix} & = \begin{bmatrix} -1 & 0 \\ 0 &-1  \end{bmatrix} \begin{bmatrix} x_1(t) \\ x_2(t) \end{bmatrix} + \begin{bmatrix} 1 \\ 2 \end{bmatrix}u(t),\\
		y &= x_1^2 + x_2^2 = \begin{bmatrix} x_1(t)&x_2(t) \end{bmatrix} \begin{bmatrix} 1 & 0 \\ 0 & 1 \end{bmatrix}\begin{bmatrix} x_1(t) \\ x_2(t) \end{bmatrix}.
		\end{align}
	\end{subequations}
	It can be easily seen that the controllability Gramian $P$ for the system \eqref{eq:2by2example} is 
	\begin{equation}
	P = \begin{bmatrix} 0.5 & 1 \\1 &2  \end{bmatrix}.
	\end{equation}
	If the method, proposed in \cite{Beeumen.2012}, is employed, then we need for solve the following Lyapunov equation for the observability Gramian:
	\begin{equation}\label{eqn:example_obs_old}
	A^T Q + QA + \begin{bmatrix} 1 & 0 \\ 0 & 1 \end{bmatrix} = 0.
	\end{equation}
	On the other hand, if we aim at employing the proposed method, \Cref{alg:SPBT_Algorithm}, then we need to solve the following equation for the observability Gramian
	\begin{align}\label{eqn:example_obs_new}
	A^T Q + QA &+ MPM = 0,
	\end{align}
	where $MPM = \begin{bmatrix} 1 & 0 \\ 0 & 1 \end{bmatrix} \begin{bmatrix} 0.5 & 1 \\1 &2  \end{bmatrix} \begin{bmatrix} 1 & 0 \\ 0 & 1 \end{bmatrix} =\dfrac{1}{2} \begin{bmatrix} 1 \\ 2  \end{bmatrix} \begin{bmatrix} 1&2   \end{bmatrix}^T$.
	Now, note that the right-hand sides of \cref{eqn:example_obs_old,eqn:example_obs_new} are of ranks $2$ and $1$, respectively. As we know from the low-rank solvers for Lyapunov equation that higher the rank of the right-hand of a Lyapunov equations, more it is expensive to determine a low-rank solution. Hence, the new proposed observability Gramian most likely be comparatively computationally cheaper. This effect in terms of computational can be seen even more when the matrix $M$ is not a positive-semi definite. 
\end{itemize}

In the next section, we prove that the proposed method preserves stability, and it possesses a guaranteed error bound,  which can be given as a function of the neglected singular values.

\section{Stability Preservation and Guaranteed Output Error Bounds}\label{sec:error_bounds} In this section, we derive some theoretical results for the BT method proposed in \Cref{subsec:SPBT}. Firstly, as for the case of \LDS~systems \cite{morPerS82} and bilinear systems \cite{morBenDRetal16}, we show that the procedure preserves stability under weak assumptions.  Secondly, we derive error bounds in the output error between the original system and the reduced-order system. For BT of  \LDS~systems,  error bounds are available, relating the $\cH_{\infty}$ norm \cite{morEnn84,morGlo84},  the $\cH_{2}$ norm \cite[Thm. 7.10]{morAnt05} and \cite{chahlaoui2012posteriori}. 
Here, we generalize the concept of $\cH_2$ norm for \LDSQO~systems, enabling us to develop a time-domain error bound with respect to the $L_{\infty}$ norm.

\subsection{Stability preservation}
It is worth noting that the classical BT for \LDS~systems produces a stable reduced-order system, see \cite[Theorem 3.2]{morPerS82}. In what follows, we provide an equivalent result for the proposed balancing method to \LDSQO~systems.  Its proof is inspired by the original one from \cite{morPerS82} and extended to \LDSQO~systems.

\begin{theorem}[Stability preservation]\label{theo:StabilCond} Suppose $H = (A, B, M)$ is a stable balanced \LDSQO~system and  $\hH:=(A_{11}, B_1, M_{11})$ is a reduced system obtained by the proposed method (\Cref{alg:SPBT_Algorithm}). Then, $A_{11}$ is also asymptotically stable  if $\Lambda(\Sigma_1)\cap \Lambda(\Sigma_2) = \emptyset$. 
\end{theorem}

\begin{proof}
	The proof is given in \Cref{apend:ProofStabil}.
\end{proof}

From \Cref{theo:StabilCond}, if $\sigma_{r} >\sigma_{r+1}$, the reduced-order system will be also stable. In practice, this condition is generally satisfied and stability is preserved.
\begin{remark} For \LDS~systems,  if $\sigma_{r} >\sigma_{r+1}$, the reduced system is  minimal. On contradictory, that is not true for \LDSQO~systems. To illustrate, let us consider the following \LDSQO~system:
	\[ A =  \begin{bmatrix}
	-1/4 & -1/3 \\ 
	-1/3 & -3/2
	\end{bmatrix}, 
	\quad B = \begin{bmatrix}
	1 & 0 \\ 
	1 & \sqrt{2}
	\end{bmatrix}, \quad \text{and} \quad M = \begin{bmatrix}
	0 & 1 \\ 1 & 1
	\end{bmatrix}. \]
	The Gramians related to the system are $P = Q = \begin{bmatrix}
	2 & 0 \\ 0 & 1
	\end{bmatrix}$. Hence, the reduce-order system associated to the highest singular value is $\hat{A} = -\dfrac{1}{4}, \hat{B} = \begin{bmatrix}
	1 & 0
	\end{bmatrix}$ and $\hat{M} = 0$, which is clearly not minimal but a zero system. 
\end{remark}
In what follows, we derive error bounds for the output approximation.

\subsection{First error bound expression} In this section, we begin by defining the notions of $\cH_2$ norm and inner product for \LDSQO~systems. Based on this, we develop an error bound for the output approximation. First, recall that the system output is given by 
\begin{align*} 
y(t) &= x^T(t)Mx(t) = (x^T(t)\otimes x^T(t))\vecop{M}\\ 
&=\vecop{M}^T(x(t)\otimes x(t)) \\
&=\int_0^t\int_0^t  \vecop{M}^T\left(e^{A\sigma_1}Bu(t-\sigma_1)\otimes e^{A\sigma_2}Bu(t-\sigma_2)\right)\diff\sigma_1\diff\sigma_2  \\
&=\int_0^t\int_0^t  \vecop{M}^T\left(e^{A\sigma_1}B\otimes e^{A\sigma_2}B\right)\left(u(t-\sigma_1)\otimes  u(t-\sigma_2)\right)\diff\sigma_1\diff\sigma_2 \\
& =\int_0^t\int_0^t  h(\sigma_1,\sigma_2)\left(u(t-\sigma_1)\otimes  u(t-\sigma_2)\right)\diff\sigma_1\diff\sigma_2, 
\end{align*}
where 
\begin{equation}\label{eq:QO_Kernel}
h(\sigma_1,\sigma_2) =  \vecop{M}^T\left(e^{A \sigma_1}B\otimes e^{A \sigma_2}B\right) = B^Te^{A^T \sigma_1}M e^{A \sigma_2}B.
\end{equation}

Hence, an \LDSQO~system can be rewritten as a 2-D convolution whose kernel is given by~\eqref{eq:QO_Kernel}. 
Consequently, we have 
\begin{align*}
|y(t)| &\leq \int_0^t\int_0^t \left| h(\sigma_1,\sigma_2)\left(u(t-\sigma_1)\otimes  u(t-\sigma_2)\right)\right|\diff\sigma_1\diff\sigma_2\\
&\leq \int_0^{t}\int_0^{t} \left\|  h(\sigma_1,\sigma_2)\diff\sigma_1\diff\sigma_2\right\|_F\left\|\left(u(t-\sigma_1)\otimes  u(t-\sigma_2)\right)\right\|_2\diff\sigma_1\diff\sigma_2 \\
&\leq \left(\int_0^{t}\int_0^{t} \left\|  h(\sigma_1,\sigma_2)\right\|_F^2\diff\sigma_1\diff\sigma_2\right)^{\frac{1}{2}} \left(\int_0^{t}\int_0^{t} \left\|u(\sigma_1)\otimes  u(\sigma_2)\right\|_2^2\diff\sigma_1\diff\sigma_2\right)^{\frac{1}{2}}, 
\end{align*}
where the last relation is followed by the Cauchy-Schwartz inequality. Additionally, given an original system with the kernel $h(\sigma_1,\sigma_2) = B^Te^{A^T \sigma_1}M e^{A \sigma_2}B$ and a reduced-order system with a  kernel $\hh(\sigma_1,\sigma_2) = \hB^Te^{\hA^T \sigma_1}\hM e^{\hA \sigma_2}\hB$, the output error can be bounded as
\begin{equation}\label{eq:EB}
\|y-\hy\|_{L_{\infty}}\leq \left(\int_0^{\infty}\int_0^{\infty} \left\|  h(\sigma_1,\sigma_2) - \hh(\sigma_1,\sigma_2)\right\|_F^2\diff\sigma_1\diff\sigma_2 \right)^{\frac{1}{2}}\|u\otimes u\|_{L_2}, 
\end{equation}
where $\|y\|_{L_{\infty}} = \max_{t\geq 0}|y(t)|$. This leads us to the following definition.
\begin{definition}[$\cH_2$ norm and inner product for \LDSQO~systems]\label{def:h2innernorm}
	Let $H = (A, B, M)$ and $\hH = (\hA, \hB, \hM)$ be stable \LDSQO~systems. Then, the $\cH_2$ norm of $H$ is defined as 
	\begin{equation}\label{eq:H2norm1}
	\|H\|_{\cH_2} = \sqrt{\left(\int_0^{\infty}\int_0^{\infty} \left\|  h(\sigma_1,\sigma_2)\right\|_F^2\diff\sigma_1\diff\sigma_2\right)^{\frac{1}{2}}}, 
	\end{equation}
	and the $\cH_2$ inner product is defined as
	\begin{equation}\label{eq:H2norm2}
	\langle H, \hH \rangle_{\cH_2}  = \int_0^{\infty}\int_0^{\infty} \trace{h(\sigma_1,\sigma_2)\hh(\sigma_1,\sigma_2)^T}\diff\sigma_1\diff\sigma_2 , 
	\end{equation}
	where $h(\sigma_1,\sigma_2) = B^Te^{A^T \sigma_1}M e^{A \sigma_2}B$ and $\hh(\sigma_1,\sigma_2) = \hB^Te^{\hA^T \sigma_1}\hM e^{\hA \sigma_2}\hB$.
\end{definition}

Also, the $\cH_2$ norm and inner-product can be characterized by the Sylvester equation, which is provided in the following proposition. 

\begin{proposition}[$\cH_2$ norm and inner product for \LDSQO~systems]\label{prop:H2normCharacterization}
	Let $H = (A, B, M)$ and $\hH = (\hA, \hB, \hM)$ be stable \LDSQO~systems of order $n$ and $r$, respectively. Then, the $\cH_2$ inner product between the two systems can be characterized as
	\begin{equation}
	\langle H, \hH \rangle_{\cH_2} = \trace{B^TZ\hB}, 
	\end{equation}
	where $Z\in \R^{n \times r}$ is a unique solution of the following Sylvester equation:
	\begin{equation}\label{eq:SylvesterObserv}
	A^TZ+Z\hA +MX\hM = 0
	\end{equation}
	in which  $X\in \R^{n \times r}$ is also a unique solution of the Sylvester equation, given by
	\begin{equation}\label{eq:SylvesterReach}
	AX+X\hA^T +B\hB^T = 0.
	\end{equation}
	Moreover, the  $\cH_2$-norm of $H$ can readily be characterized as
	\begin{equation*}
	\|H\|_{\cH_2}= \sqrt{\trace{B^TQB}}, 
	\end{equation*}
	where $Q\in \R^{n \times n}$ is the observability Gramian.
\end{proposition}
\begin{proof}
	Since both $H$ and $\hH$ are stable systems, the Sylvester equations \eqref{eq:SylvesterReach} and  \eqref{eq:SylvesterObserv} have unique solutions which can also be given by
	\[X= \int_{0}^{\infty} e^{A\sigma_1}B\hB e^{\hA^T\sigma_1} d\sigma_1\quad \text{and}\quad  Z= \int_{0}^{\infty} e^{A^T\sigma_2}MX\hM e^{\hA\sigma_2} d\sigma_2.  \]
	By inserting the expression of $X$ into $Z$ and following by multiplying with $B^T$ and $\hB$ from the left and right-hand sides, respectively, we have the desired result.  
\end{proof}

Finally, by assembling  the results together and notions above, the following a posteriori error bound holds.
\begin{theorem}[A posteriori error bound]\label{theo:AposterioriEB} Let $H = (A, B, M)$ and $\hH = (\hA, \hB, \hM)$ be stable \LDSQO~systems and suppose that $y$ and $\hat{y}$ are their respective outputs, subject to the same input $u$. Then,
	\begin{equation}\label{eq:AprioriEB}
	\begin{array}{rl}
	\|y-\hat{y}\|_{L_{\infty}} &\leq \|H-\hH\|_{\cH_2}\|u\otimes u\|_{L_2}
	\end{array}
	\end{equation}
	where
	\begin{equation}\label{eq:SylvesterEB} 
	\|H-\hH\|_{\cH_2} =  \sqrt{\trace{B^TQB +\hB\hQ\hB -2B^TZ\hB}}, \end{equation}
	where $Q$ and $\hat{Q}$ are the observability Gramians for the systems $H$ and $\hH$, respectively, and $Z$ is the cross-Gramian, which is the solution of \cref{eq:SylvesterObserv}.
\end{theorem}
\begin{proof}
	\Cref{eq:AprioriEB} corresponds to \cref{eq:EB}  with the norm notation. Additionally, we make use of the triangular inequality, i.e., $\|H-\hat{H}\|^2_{\cH_2} = \|H\|_{\cH_2}^2+\|\hH\|_{\cH_2}^2 -2	\langle H, \hH \rangle_{\cH_2},$ and  \cref{prop:H2normCharacterization}, resulting into  \cref{eq:SylvesterEB}. 
\end{proof}

\Cref{theo:AposterioriEB} provides a bound for the output error between the original system and the reduced-order one. This bound relies on the $\cH_2$ norm of the error between the $H$ and $\hH$. Hence, if they are close with respect to the $\cH_2$ norm, they will both provide outputs that are also close. Notice that the results presented in this subsection do not assume that a reduced-order system is obtained by the proposed  BT procedure. Indeed, \Cref{theo:AposterioriEB} holds for any reduced-order system $\hH = (\hA, \hB, \hM)$, provided that the matrix $\hA$ is stable.  In the next subsection, we will show for the proposed BT how this error bound relates to the singular values. 
\subsection{Error bound and singular values} In this section, we will show how the error bound given in \cref{theo:AposterioriEB} is related to the  singular values in the proposed BT procedure. To that aim, suppose that $H = (A, B, M)$ is an $n$-order balanced system, whose  Gramians $P = Q= \Sigma = \diag(\sigma_1, \dots, \sigma_n)$. Furthermore, we consider that the system matrices are partitioned as~\eqref{eq:PartionedMatrices}. 
Next, we present how error can be bounded in terms of the  singular values. This result was mainly inspired by the error bounds for \LDS~systems regarding the $\cH_2$ and time-limited $\cH_2$ norms, see \cite[Thm. 7.10]{morAnt05} and \cite{chahlaoui2012posteriori, redmann2018output, Pontes2019numerical}. 
\begin{theorem}[BT error bound as a function of the singular values]\label{theo:EBhsv}
	The $\cH_2$ norm of the error system is given by	
	\begin{equation}
	\|\cE\|_{\cH_2}^2 = \trace{\left(B_2B_2^T+  2Z_2A_{12}+2M_{:2}^TXM_{12}+2X_2A_{21}^T\right)\Sigma_2} + \trace{B_1^T\right(\hat{Q}-\Sigma_1\left)B_1}, 
	\end{equation}
	where $M_{:2}^T = \begin{bmatrix}
	M_{12}^T & M_{22}
	\end{bmatrix}$. 
	where $M_{:2}^T = \begin{bmatrix}
	M_{12}^T & M_{22}
	\end{bmatrix}$.

\end{theorem} 
\begin{proof}
	The proof is given in Appendix B.
\end{proof}
\cref{theo:EBhsv} shows how the error bound for BT relates with the neglected  singular values $\Sigma_2$ and the preserved singular values $\Sigma_1$. The following corollary removes the dependency of $\Sigma_1$ in the error bound.
\begin{corollary}[BT error bound depending on neglected singular values]\label{cor:EBnegHSV}  The following error bound holds
	\begin{equation}
	\|\cE\|_{\cH_2}^2 \leq \trace{\left(B_2B_2^T+  2Z_2A_{12}+2M_{:2}^TXM_{12}+2X_2A_{21}^T\right)\Sigma_2}.
	\end{equation}
\end{corollary}
\begin{proof}
	To prove the above result, we show that $\trace{B_1^T\right(\hat{Q}-\Sigma_1\left)B_1} \leq 0$. This result follows because  $\hQ-\Sigma_1$ is a negative semi-definite matrix. Indeed, $\cD := \hQ-\Sigma_1$ satisfies the following Lyapunov equation:
	\[A_{11}^T\cD +\cD A_{11} -M_{12}\Sigma_2M_{12}^T = 0. \]
	Since $A_{11}$ is Hurwitz and $-M_{12}\Sigma_2M_{12}^T$ is a negative semi-definite matrix, $\cD$ is negative semi-definite. This proves the result.
\end{proof}
\cref{cor:EBnegHSV} shows that the $\|\cE\|_{\cH_2}$ can be bounded only by the terms from $\Sigma_2$. Furthermore, the term $\trace{B_1^T\right(\hat{Q}-\Sigma_1\left)B_1}$ can be written as function of $\Sigma_2$, since 
\[ \cD = -\int_0^{\infty} e^{A_{11}^T \sigma}M_{12}\Sigma_2M_{12}^Te^{A_{11}\sigma} d\sigma \] and
$\trace{B_1^T\right(\hat{Q}-\Sigma_1\left)B_1} = \trace{B_1^T\cD B_1}$.


\subsection{Small scale example}

Now we illustrate the obtained results by applying the proposed method to a small-scale system and by computing the error bounds from \cref{theo:AposterioriEB} (or \Cref{theo:EBhsv}) and  \cref{cor:EBnegHSV}. 

\begin{figure}[H]
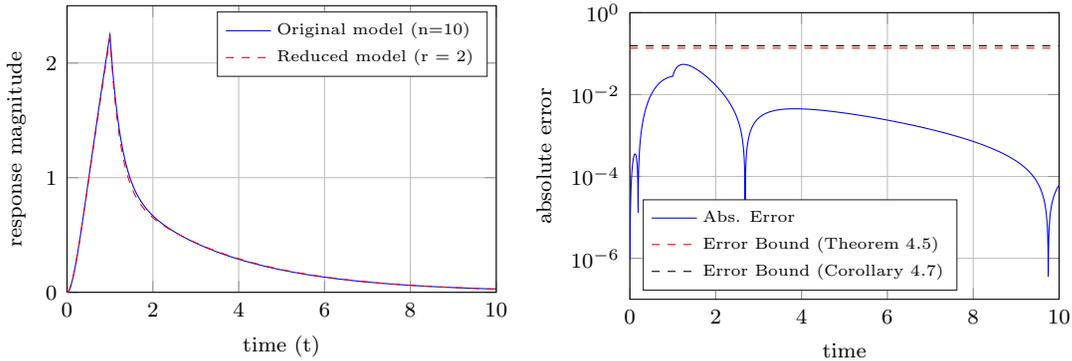

	\begin{center}
		\input{TimeDomain_SmallEx10to6.tex}
		\input{AbsErr_SmallEx10to6.tex}
		\caption{Small sclae example: Output errors $|y(t)-\hy(t)|$ and error bounds from \cref{theo:AposterioriEB} and \cref{cor:EBnegHSV} for  small scale example to order $n=10$, reduced order $r =2$.}
		\label{fig:SmallArticleEx}
	\end{center}
\end{figure}

We consider a random stable single-input single-output (SISO) system of order $n = 10$, generated by using the command \textsf{rss}  in \matlab\xspace  with \texttt{seed} $0$. The reduced-order system of order $r=2$ is computed using the proposed method, and the error bounds are computed using \matlab~ direct solver (command \textsf{lyap}). We simulate the time domain response of  the corresponding the original and reduced-order systems, using an input  $u(t) = e^{-\frac{1}{4}t}$ for $t\geq 0$.  In this case, $\|u\otimes u\|_{L_2} = \|u^2\|_{L_2} = 1$. The results of the absolute errors are depicted in \Cref{fig:SmallArticleEx}, as well as the  bounds from  \cref{theo:AposterioriEB} and \cref{cor:EBnegHSV}.

As expected, the proposed method produces an stable reduced-order system. Additionally, by inspecting the time-domain error, we observe that it satisfies the proposed error bounds.

\section{Numerical Experiments}\label{sec:experiments}
In this section, we test the efficiency of the proposed methodology, described in \Cref{alg:SPBT_Algorithm} (denoted by \texttt{SPBT}) for \LDSQO~systems by means of two numerical examples and compare with the existing methods as discussed in \Cref{alg:LTBT_Algorithm} (denoted by \texttt{LTBT}), and \Cref{alg:SPBT_Algorithm} (denoted by \texttt{QBTBT}).  All the simulations are done using \matlab~R2017a (64-bit) on a machine with \intel\coreifive-6600 processors with 3.3-GHz clock frequency, 8 GB RAM and Windows 8 operating system.  To solve the Lyapunov equations,  we use the ADI-solvers provided in the {M-M.E.S.S.} toolbox~\cite{saak_2016}.

\subsection{Clamped beam}
As a first example, we discuss a clamped beam model. It is widely used as one of benchmark MOR problems, see e.g.,~\cite{slicot_beam}.  A detailed description of the dynamics can be found in the mentioned reference; therefore, we omit it in the interest of brevity.  However, we consider a variant of the example by only modifying the output equation while keeping the differential equation for the state vector $x(t)$ same. For the output equation, we define a diagonal matrix $M$ such that the output $x(t)^TMx(t)$ is a weighted sum of the squares of $100$ randomly components of the state. The random components are selected by setting $\texttt{seed}=1$. All the weights are the same and the sum of the weights is equals to 1. The order of the system is $348$.

\begin{figure}[!tb]
	\centering
	\input{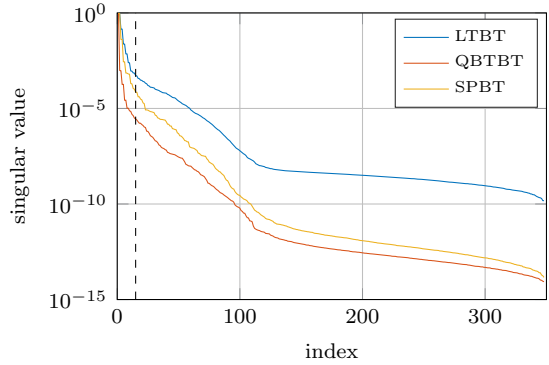} 
	\caption{Clamped Beam: Comparison of the normalized singular values.
	}
	\label{fig:sim_det_singularvalues}
\end{figure}

Next, we aim at employing \texttt{SPBT}, \texttt{LTBT} and \texttt{QBTBT} to compute reduced-order systems. For this,  we plot the decay of the normalized singular values, obtained by the all three considered methods in \Cref{fig:sim_det_singularvalues}. Next, we determine reduced systems of order $r = 15$. Also,  for this example, the $\cH_2$ norm of the original system is  $\|H\|_{\cH_2} \approx 5.12\cdot 10^3 $, the norm of the system error $\|H-\hH\|_{\cH_2} \approx 8.88$ and the relative norm of the error system is $\|H-\hH\|_{\cH_2}/\|H\|_{\cH_2}  \approx 2\cdot 10^{-3}$. 
To compare the quality of these approximations, we simulate the original system and the obtained reduced-order systems using a sinusoidal excitation $u_1(t) = \sin\left(\frac{2 \pi t}{10}\right)+1$  and an exponentially damped quadratic excitation $u_2(t) =  \left(e^{-\frac{t}{5}} \right) t^2$. In \Cref{fig:sim_det_response}, we plot the transient responses of the original and reduced-order systems, and also the absolute errors between the outputs of the original and reduced-order systems for the considered input signals.   These figures show that \texttt{SPBT} outperforms \texttt{QBTBT}; on the other hand, we see that \texttt{SPBT} is very competitive with \texttt{LTBT} for this example. Moreover, \texttt{SPBT} can be shown to be more efficient than \texttt{LTBT} as it includes the factor $MPM$ in the Lyapunov equation for the observability Gramian which can reduce the solution search space for the observability Gramian to the space of controllable states. 

Furthermore, we compute the error bound for $\|y(t) -\hy(t)\|_{L_\infty}$ using \Cref{theo:AposterioriEB}. Since  this error bound assume that the input is in $L_2$, we are only able to compute it for the case of exponentially damped quadratic excitation. Hence, for this case, we  compute the error bound and show in \Cref{fig:clampedBeamu2}. On the other hand, the quantity cannot be estimated for the sinusoidal input since it is unbounded in $L_2$.

\begin{figure}[!tb]
	\begin{subfigure}[!htb]{\textwidth}
		\centering
		\input{Fig_modelDeterministic_5.tex} 
		\input{Fig_modelDeterministic_7.tex} 
		\caption{Comparison of systems' responses for the input $u_1$.}
	\end{subfigure}
	\begin{subfigure}[!htb]{\textwidth}
		\centering
		\input{Fig_modelDeterministic_6.tex} 
		\input{Fig_modelDeterministic_9.tex} 
		\caption{Comparison of systems' responses for the input $u_2$.}
		\label{fig:clampedBeamu2}
	\end{subfigure}
	\label{fig:sim_det_response}
\end{figure}

\subsection{Steel Profile}
As a second example, we consider a semi-discretized heat transfer problem for the optimal cooling of a steel profile, whose detailed description of the dynamics can be found in~\cite{morwiki_steel}. Typically, the output function is given as $y(t) = Cx(t)$. However, we modify the output function and our interest of quantity is the $2-$norm of the $y(t)$, i.e., $x(t)C^TCx(t) =: x(t)^TMx(t)$. Moreover, the system is converted into a stochastic system by augmenting the input matrix $B$ with an additional column $b_w$ for the the noise input $w$. All the entries of the vector $b_w$ are the same and the values is set as the maximum entry of the matrix $B$. This implies that when a noise signal is applied to the system, it affects each state with the same way. 
The order of the system is $n = 1357$ and the system can be represented by 
\begin{subequations}
	\label{eqn:LDS_QO_Stochastis_SystemEquations}
	\begin{align}
	\dot{x}(t) &=Ax(t) + \begin{bmatrix}B& b_w \end{bmatrix} \begin{bmatrix}u(t)\\w(t)\end{bmatrix}, \\
	y(t) &={x(t)}^{T}Mx(t),
	\end{align}
\end{subequations}

\begin{figure}[tb]
	\centering
%
\definecolor{mycolor1}{rgb}{0.00000,0.44700,0.74100}%
\definecolor{mycolor2}{rgb}{0.85000,0.32500,0.09800}%
\definecolor{mycolor3}{rgb}{0.92900,0.69400,0.12500}%
\begin{tikzpicture}

\begin{axis}[%
width=0.951\figurewidth,
height=\figureheight,
at={(0\figurewidth,0\figureheight)},
scale only axis,
xmin=0,
xmax=72,
xlabel style={font=\color{white!15!black}},
xlabel={index},
ymode=log,
ymin=1e-20,
ymax=1,
yminorticks=true,
ylabel style={font=\color{white!15!black}},
ylabel={singular value},
axis background/.style={fill=white},
title style={font=\bfseries},
title={Normalized Singular Value Decay Plot},
xmajorgrids,
ymajorgrids,
yminorgrids,
legend style={legend cell align=left, align=left, draw=white!15!black},
title ={},ticklabel style={font=\scriptsize},xlabel style={font=\scriptsize},ylabel style={font=\scriptsize},legend style={font=\tiny},
]
\addplot [color=mycolor1]
  table[row sep=crcr]{%
1	1\\
2	0.356245057100573\\
3	0.200990708564166\\
4	0.149372780603414\\
5	0.133388999474348\\
6	0.0852369296257154\\
7	0.0252141863308938\\
8	0.0187678696368757\\
9	0.0163545263562626\\
10	0.00540985644189399\\
11	0.00201294890298664\\
12	0.00101306553521416\\
13	0.000240927163589824\\
14	1.69147789276567e-05\\
15	2.2543101546942e-06\\
16	1.43321812958961e-07\\
17	1.40977019298293e-08\\
18	2.83992453660295e-09\\
19	2.49440165192622e-09\\
20	9.29604887560503e-10\\
21	8.40166455049498e-10\\
22	6.44385813582027e-10\\
23	5.93571964263242e-10\\
25	3.06295361516904e-10\\
26	2.19214172073004e-10\\
27	1.59783944786093e-10\\
28	1.38678462786011e-10\\
29	1.16295132754723e-10\\
30	8.21229261082342e-11\\
32	7.29439512201144e-11\\
33	3.78303709171603e-11\\
34	2.74344558758294e-11\\
35	1.77179086972079e-11\\
36	1.43178119741508e-11\\
37	7.50791952766296e-12\\
38	7.22973178179398e-12\\
39	7.11381825823675e-12\\
40	5.64769292052196e-12\\
41	4.20322359279932e-12\\
43	1.96676102621713e-12\\
44	1.83593919311773e-12\\
45	1.06919982277511e-12\\
46	8.59183887586819e-13\\
47	4.23466564549483e-13\\
48	3.68053734366821e-13\\
49	2.36239076689225e-13\\
50	2.2633903284443e-13\\
51	2.0594769912604e-13\\
52	1.95384816416235e-13\\
53	1.18582293084621e-13\\
54	9.43610229477134e-14\\
55	6.92663142756559e-14\\
56	4.28452937648759e-14\\
57	3.20806687212595e-14\\
58	2.77365325681669e-14\\
59	2.25413713949863e-14\\
60	1.9102924143647e-14\\
61	1.56760528033623e-14\\
62	1.45898466583384e-14\\
63	8.43755049577486e-15\\
64	3.07559581956264e-15\\
65	1.59141128900201e-15\\
66	9.50421766259315e-16\\
67	8.52852271524319e-16\\
68	7.90484221120721e-16\\
69	6.48783915999459e-16\\
70	6.09331153188688e-16\\
71	1.86862943995174e-16\\
72	7.0816598798944e-17\\
};
\addlegendentry{LBT}

\addplot [color=mycolor2]
  table[row sep=crcr]{%
1	1\\
2	3.97053477273602e-05\\
3	1.34399568542317e-05\\
4	5.82832068848388e-06\\
5	4.47772676114676e-06\\
6	2.45195411408742e-06\\
7	9.48082530234625e-07\\
8	7.625405723083e-07\\
9	4.09575557774324e-07\\
10	2.07186941134256e-07\\
11	3.59683466709862e-08\\
12	1.96279428514416e-08\\
13	1.39543163164103e-08\\
14	4.91922720766093e-09\\
15	4.30929293933261e-10\\
16	5.68947551319013e-11\\
17	1.27225818594806e-12\\
18	1.30386070226182e-13\\
19	4.31601510709593e-15\\
20	9.99200722162635e-17\\
72	9.99200722162635e-17\\
};
\addlegendentry{QBBT}

\addplot [color=mycolor3]
  table[row sep=crcr]{%
1	1\\
2	0.162539327936651\\
3	0.108967553252216\\
4	0.0644068004909917\\
5	0.0320619386847544\\
6	0.0213127690733785\\
7	0.0165025878911979\\
8	0.0065632280618668\\
9	0.00250513998581142\\
10	0.000530014734602222\\
11	0.000296035846594267\\
12	0.000186141297884745\\
13	6.02911298802092e-05\\
14	4.82781835950045e-06\\
15	8.46782471289691e-07\\
16	1.54616257819369e-08\\
17	1.90131760322038e-09\\
18	9.61681992054811e-17\\
19	7.42066062800464e-17\\
71	7.42066062800464e-17\\
72	1.49073513741117e-17\\
};
\addlegendentry{SPBT}

\addplot [color=black, dashed, forget plot]
  table[row sep=crcr]{%
15	1e-20\\
15	1\\
};
\end{axis}
\end{tikzpicture}%
	\caption{Steel Profile: Comparison of the normalized singular values.
	}
	\label{fig:sim_stc_singularvalues}
\end{figure}
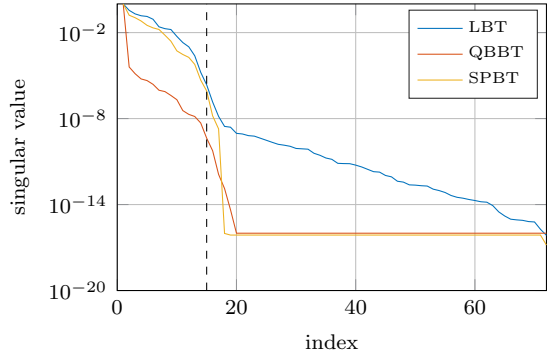

As can be seen, all the theory in the paper is developed for deterministic systems, thus all the results might not readily be extended to stochastic systems. 
Nevertheless, \Cref{thm:ContObser_PQ}, identifying of uncontrollable and unobservable subspaces, can be proven to be held even for the stochastic case. Hence, we blindly  employ \texttt{SPBT}, \texttt{LTBT} and \texttt{QBTBT}, and aim at constructing order model models. To that end, in \Cref{fig:sim_stc_singularvalues}, we  first show the decay of the  normalized singular values, where we make a similar observation as in the previous example. Next, we determine reduced-order systems of order $r=15$ using the considered methods.

To compare the quality of these approximations, we simulate the original and the reduced systems obtained by each method using setting each input component of the input vector $u(t)$ to a sinusoidal excitation $u_1(t) = 20 \sin \Big(\frac{2 \pi t}{10}\Big)+1$ and an exponentially damped quadratic excitation  $u_2(t) =  3 \big(e ^{-\frac{t}{5}} \big) t^2$. For a noisy input, a white Gaussian noise sequence with power $1\text{dB}$ and a magnitude scaling factor is set to $40$. A total of 50 simulations are  performed to get an approximation of the expected response. The systems are simulated using a $4$th order Runge-Kutta method (\texttt{RK4}) based solver. 

\begin{figure}[tb]
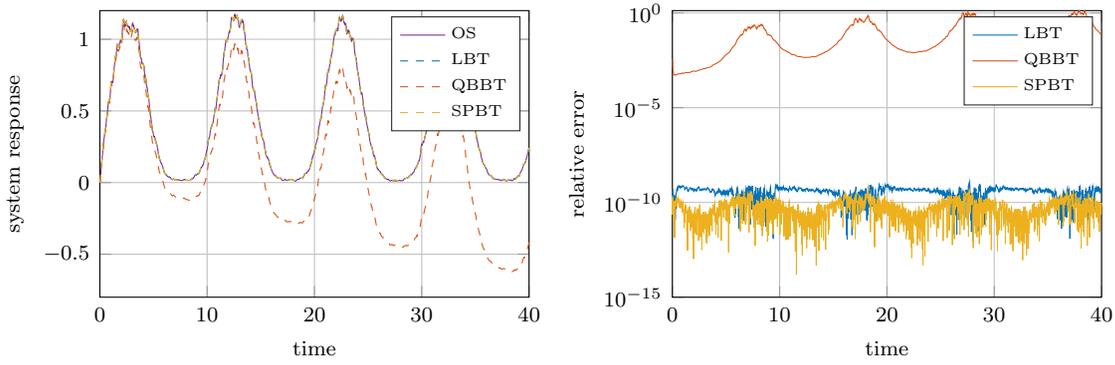
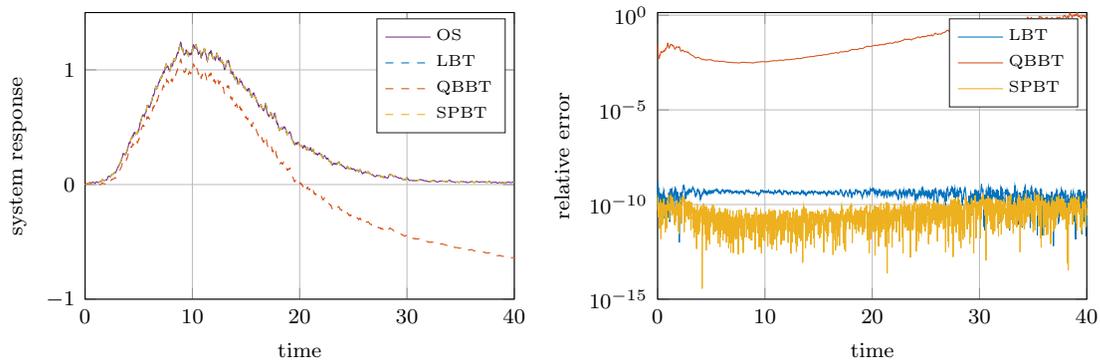

	\begin{subfigure}[!htb]{\textwidth}
		\centering
		\input{Fig_modelStochastic_6.tex} 
		\input{Fig_modelStochastic_9.tex} 
		\caption{Comparison of systems' responses to for the input $u_1$.}
	\end{subfigure}
	\begin{subfigure}[!htb]{\textwidth}
		\centering
		\input{Fig_modelStochastic_7.tex} 
		\input{Fig_modelStochastic_11.tex} 
		\caption{Comparison of systems' responses to for the input $u_2$.}
	\end{subfigure}
	\caption{Steel Profile: Comparison of the transient responses of the original and order model models obtained by various methods for two arbitrary control inputs and Gaussian noise.}
	\label{fig:sim_stc_response}
\end{figure}

\Cref{fig:sim_stc_response} presents the transient responses and response errors of the output for these input signals, which \texttt{QBTBT} fails to capture the dynamics of the system; on the other hand, we observe that the proposed \texttt{SPBT} is marginally better as compared to \texttt{LTBT} for this example. 

\section{Conclusions}\label{sec:conclusion}
In this paper, we have studied a balanced truncation method for linear dynamical systems with a quadratic output function. For this, we have proposed a new pair of algebraic Gramians. In particular, the observability Gramian has been introduced by making use of the corresponding adjoint system. Consequently, we have shown that these Gramians encode controllability and observability of the system and have studied the connection between the proposed Gramians and energy functionals. This has allowed us to proposed an algorithm to determine reduced-order systems by truncating unimportant states for the input-output dynamics.  We have also discussed advantages of the proposed methods over the existing methods in the literature. Finally, based on  $\mathcal H_2$ energy considerations, we have derived error bounds, depending on the neglected singular values.  Furthermore, we have shown an efficiency of the proposed methods by means of a couple of numerical examples. 
%

\section*{Acknowledgements}
We would like to thank Adil Ahsan for helping us with the \matlab~implementations at Max Planck Institute for Dynamics of Complex Technical Systems, Magdeburg, Germany.

\bibliographystyle{siamplain}
\bibliography{mor,csc}

\appendix
\section{Proof of \cref{theo:StabilCond}}\label{apend:ProofStabil}
For the balanced realization, \eqref{eq:RedRechGram_Part} holds, i.e., 
\begin{equation}\label{eq:RedReachGramTheo}  A_{11}\Sigma_1 +\Sigma_1A_{11}^T +B_1B_1^T = 0.
\end{equation}
This implies, from $\Sigma_1>0$ and  \cite[Lemma 3.1]{morPerS82}, that the eigenvalues of $A_{11}$ should have real part less or equal to zero\footnote{This result can also be seen as a consequence of the Inertia Theorem.}.  However, we still need to show that $A_{11}$ has no eigenvalues on the imaginary axis. 

Suppose, by contradiction, that  $\nu = i \omega$, $\omega \in \R$ is a  eigenvalue of $A_{11}$ and $U$ is a basis of the kernel of $(\nu I - A_{11})$. Hence,
\[ A_{11}U = \nu U\quad \text{and} \quad U^*A_{11}^T = \nu^* U^*,  \]  
where $\nu^*$ denotes the complex conjugate of $\nu$ and $U^*$ denotes the Hermitian transpose of $U$. From $\eqref{eq:RedObservGram_Part}$, we have
\begin{equation}\label{eq:RedObservGramTheo}
A_{11}^T\Sigma_1 +\Sigma_1A_{11} +M_{11}\Sigma_1 M_{11} +M_{12}\Sigma_2 M_{12}^T=0. 
\end{equation}
Multiplying \eqref{eq:RedObservGramTheo} from the right and left-hand sides with $U$ and $U^*$, receptively, yields
\[U^*M_{11}\Sigma_1M_{11}U  + U^*M_{12}\Sigma_2M_{12}^TU =0.\]
Since $\Sigma_1, \Sigma_2>0$, this yields 
\[M_{11}U = 0 \quad \text{and} \quad M_{12}^TU = 0.\]  
Now, by multiplying \eqref{eq:RedObservGramTheo} by $U$ from the right-hand side, we obtain 
\[ (A_{11}^T+\nu I_r)\Sigma_1U = 0.\] 
Furthermore, by multiplying \eqref{eq:RedReachGramTheo} from the right and left-hand sides by $\Sigma_1 U$ and $U^*\Sigma_1$, we obtain 
\[U^*\Sigma_1BB^T\Sigma_1U = 0\implies B^T\Sigma_1U = 0.\] 
Again, by multiplying  \eqref{eq:RedReachGramTheo} from the right-hand side with $\Sigma_1 U$, we get \[(A_{11}-\nu I_r)\Sigma_1^2U = 0.\] 
As a result, it can be noticed that $\Sigma_1^2U$ is also a basis fo the null-space of $(\nu I -A_{11})$, thus leading to
$\Sigma_1^2U = U\overline{\Sigma}_1^2$, for some matrix $\overline{\Sigma}_1^2$. Note that the eigenvalues of $\overline{\Sigma}_1^2$ are a subset of those from $\Sigma_1^2$. From the structure of the problem, it is possible to choose $U$ such that $\overline{\Sigma}_1^2$ is diagonal, whose elements are a subset of the diagonal entries of $\Sigma_1$. 

Next, we consider the blocks $(2,1)$ from  \eqref{eq:FullReachGram} and \eqref{eq:FullObsevGram}, which are as follows:
\begin{subequations}
	\begin{align}
	A_{21}\Sigma_1 +\Sigma_2 A_{12}^T + B_2B_1^T &= 0, \label{eq:eq21}\\ 
	A_{12}^T\Sigma_1 +\Sigma_2 A_{21}+ M_{12}^T\Sigma_1M_{11} +M_{22} \Sigma_2M_{12}^T &= 0.\label{eq:eq22}	 
	\end{align}
\end{subequations}
By multiplying  \cref{eq:eq21} by $\Sigma_1U$ and  \cref{eq:eq22} by $U$ from the right-hand side, we have
\begin{align}
A_{21}\Sigma_1^2U +\Sigma_2A_{12}^T\Sigma_1U &= 0,\label{eq:eqFinal1} \\ 
A_{12}^T\Sigma_1 U +\Sigma_2 A_{21}U  &= 0.\label{eq:eqFinal2}
\end{align}
Multiplying  \eqref{eq:eqFinal2} on the left-hand side by $\Sigma_2$ and subtracting \eqref{eq:eqFinal1} yields
\[\Sigma_2^2A_{21}U = A_{21}\Sigma_1^2 U = A_{21}U\overline{\Sigma}_1^2. \]
Hence, one can write
\[\begin{bmatrix}
\overline{\Sigma}_1^2 & 0 \\ 0 &\Sigma_2^2
\end{bmatrix} \begin{bmatrix}
I \\ A_{21}U
\end{bmatrix} = \begin{bmatrix}
I \\ A_{21}U
\end{bmatrix} \overline{\Sigma}_1. \]
From the hypothesis that $\overline{\Sigma}_1^2$ and $\Sigma_2^2$ have no common eigenvalues, it follows that $A_{21}U= 0$. As a consequence, we have 
\[ \begin{bmatrix}
A_{11} & A_{12} \\ 
A_{21} & A_{22}
\end{bmatrix} \begin{bmatrix}
U \\ 0 
\end{bmatrix} = \nu \begin{bmatrix}
U \\ 0
\end{bmatrix}, \]
which contradicts the fact that the original matrix A is Hurwitz. Hence, the matrix $A_{11}$ is also Hurwitz.

\section{Proof of \cref{theo:EBhsv}}
Let us consider the matrices $X$ and $Z$, which are the solutions of the Sylvester equations \eqref{eq:SylvesterReach} and \eqref{eq:SylvesterObserv}, respectively, and are portioned in the same way, i.e.,
$Z^T = \begin{bmatrix}
Z_1^T & Z_2^T
\end{bmatrix}$ and $X^T = \begin{bmatrix}
X_1^T & X_2^T
\end{bmatrix}$, with  $Z_1, X_{1}\in \R^{r \times r}$ and $Z_2, X_2\in \R^{(n-r)\times r}$,
Hence, the $\cH_2$ norm of the error system can be given by
\begin{align} \|\cE\|_{\cH_2}^2 &= \trace{B^T\Sigma B - 2B^TZB_1 + B_1^T\hat{Q}B_1} \nonumber \\
&= \trace{B^T\Sigma B - 2B_1^TZ_1B_1-2B_2^TZ_2B_1 + B_1^T\hat{Q}B_1},
\label{eq:FirstNormDev}
\end{align}
where $\hQ$ satisfies
\[A_{11}^T\hat{Q} + \hat{Q}A_{11} +M_{11}\Sigma_1M_{11} = 0. \]
Firstly, let us analyze the term $\trace{2B_2^TZ_2B_1}$. By developing the block (2,1) of \cref{eq:FullReachGram},  we obtain
\[A_{12}\Sigma_2 +\Sigma_1A_{21}^T+B_1B_2^T= 0,\]
and, consequently, $-B_1B_2^TZ_2 = (A_{12}\Sigma_2 +\Sigma_1A_{21}^T)Z_2$. This implies
\[ \trace{-2B_2^TZ_2B_1} = \trace{-2B_1B_2^TZ_2} = \trace{2A_{12}\Sigma_2Z_2 +2\Sigma_1A_{21}^TZ_2}.\]
Substituting the above relation in \eqref{eq:FirstNormDev} yields
\begin{align*}
\|\cE\|_{\cH_2}^2 = &\trace{B_2^T\Sigma_2B_2+  2A_{12}\Sigma_2Z_2} + \trace{B_1^T\right(\hat{Q}-\Sigma_1\left)B_1} \\&+2\trace{B_1^T\Sigma_1B_1-B_1^TZ_1B_1+\Sigma_1A_{12}^TZ_2}.
\end{align*}
Now, let us have a closer look at the block $(1,1)$ of $A^TZ+ZA_{11} +MX M_{11}$, that is:
\[A^T_{11}Z_1 +A_{21}^TZ_2+Z_1A_{11}+M_{11}X_1M_{11}+M_{12}X_2M_{11} = 0. \]
Hence, we have
\[\Sigma_1A_{12}^TZ_2 = -(\Sigma_1A_{11}^TZ_1+\Sigma_1Z_1A_{11}+ Q_o),\]
where $Q_o := \Sigma_1M_{11}X_1M_{11}+\Sigma_1M_{12}X_2M_{11}$. Thus,  we obtain
\begin{align*}
&\trace{B_1^T\Sigma_1B_1-B_1^TZ_1B_1+\Sigma_1A_{12}^TZ_2} \\
&\qquad= \trace{B_1^T\Sigma_1B_1-B_1^TZ_1B_1-\Sigma_1A_{11}^TZ_1-\Sigma_1Z_1A_{11}-Q_o} \\
&\qquad=\trace{B_1^T\Sigma_1B_1-Q_0}-\trace{\underbrace{\left(B_1B_1^T+\Sigma_1A^T_{11}+A_{11}\Sigma_1\right)}_{=0}Z_1}.
\end{align*}
Having combined all, we finally have
\begin{align*}
\|\cE\|_{\cH_2}^2 = &\trace{B_2^T\Sigma_2B_2+  2A_{12}\Sigma_2Z_2} + \trace{B_1^T\right(\hat{Q}-\Sigma_1\left)B_1}+2\trace{B_1^T\Sigma_1B_1-Q_0}.
\end{align*} 
Now, we study the following term:
\[\trace{B_1^T\Sigma_1B_1-Q_0} =\trace{B_1^T\Sigma_1B_1 - \Sigma_1M_{11}X_1M_{11}-\Sigma_1M_{12}X_2M_{11}}. \]
Moreover, by analyzing the block (1,1) of \cref{eq:SylvesterReach}, we obtain
\[A_{11}X_1+X_1A_{11}^T +A_{12}X_2+B_1B_1^T = 0\]
and recall that
\[A_{11}^T\Sigma_1+\Sigma_1A_{11}  +M_{11}\Sigma_1M_{11} + M_{12}\Sigma_2M_{12}^T= 0.\]
As a consequence, 
\[\trace{X_1^T\left(M_{11}\Sigma_1M_{11} + M_{12}\Sigma_2M_{12}^T\right)} = \trace{\Sigma_1\left(A_{12}X_2+B_1B_1^T \right)}.  \]
Hence,
\[\trace{B_1^T\Sigma_1B_1 - \Sigma_1M_{11}X_1M_{11}} = \trace{M_{12}^TX_1^TM_{12}\Sigma_2-\Sigma_1A_{12}X_2}.\]
Now, the total error can be written as 
\begin{align*}
\|\cE\|_{\cH_2}^2 = &\trace{B_2^T\Sigma_2B_2+  2A_{12}\Sigma_2Z_2+2M_{12}^TX_1^TM_{12}\Sigma_2} + \trace{B_1^T\right(\hat{Q}-\Sigma_1\left)B_1}\\
&-2\trace{\Sigma_1A_{12}X_2+\Sigma_1M_{12}X_2M_{11}}.
\end{align*} 
Finally, from the block (1,2) of $A^T\Sigma+\Sigma A+M\Sigma M$, we have
\[A_{21}^T\Sigma_2+\Sigma_1A_{12}+M_{11}\Sigma_1M_{12}+M_{12}\Sigma_2M_{22} = 0,\] 
giving us 
\[ \trace{(\Sigma_1A_{12}+M_{11}\Sigma_1M_{12})X_2} = -\trace{A_{21}^T\Sigma_2X_2+M_{12}\Sigma_2M_{22}X_2}. \]
By substituting the latter equation in the error expression,  we obtain 
\begin{align*}
\|\cE\|_{\cH_2}^2 = &\trace{\left(B_2B_2^T+  2Z_2A_{12}+2
	M_{12}^TX_1M_{12}+M_{22} X_2M_{12} +2X_2A_{21}^T\right)\Sigma_2}\\ &+ \trace{B_1^T\right(\hat{Q}-\Sigma_1\left)B_1},  
\end{align*}
which proves the theorem.

%

\end{document}